\documentclass[12pt,leqno]{amsart}
\usepackage{amssymb}
\usepackage{amsmath,amssymb,color}
\numberwithin{equation}{section}

\newcommand{\weight}{e^{2s\varphi}}

\newcommand{\ep}{\varepsilon}
\newcommand{\la}{\lambda}
\newcommand{\va}{\varphi}
\newcommand{\ppp}{\partial}

\newcommand{\www}{\widetilde}

\newcommand{\R}{\mathbb{R}}

\newcommand{\ooo}{\overline}
\newcommand{\OOO}{\Omega}

\newcommand{\sumij}{\sum_{i,j=1}^d}

\newcommand{\HHHH}{H^{2,1}}
\newcommand{\DDDDD}{D(u(\cdot,T))}

\newcommand{\hhalf}{\frac{1}{2}}

\allowdisplaybreaks

\title
[]
{
Stability in 
determination of states for the mean field game equations
}

\author{$^1$ Hongyu Liu and 
$^2$ Masahiro Yamamoto
}

\thanks{
$^1$ Department of Mathematics, City University of Hong Kong, Kowloon, 
Hong Kong SAR, P.R. China. E-mail: {\tt hongyliu@cityu.edu.hk}
\\
$^2$ Graduate School of Mathematical Sciences, The University
of Tokyo, Komaba, Meguro, Tokyo 153-8914, Japan. 
E-mail: {\tt myama@ms.u-tokyo.ac.jp}
}

\date{}
\begin{document}
\maketitle

\begin{abstract}
We consider solutions satisfying the Neumann zero boundary condition and
a linearized mean field game system in $\OOO \times (0,T)$, 
where $\OOO$ is a bounded domain in $\R^d$ and $(0,T)$ is the time
interval.  We prove two kinds of 
stability results in determining the solutions.  The first is H\"older 
stability in time interval $(\ep, T)$ with arbitrarily fixed $\ep>0$ by 
data of solutions in $\OOO \times \{T\}$.
The second is the Lipschitz stability in $\OOO \times (\ep, T-\ep)$ by 
data of solutions in arbitrarily given subdomain of $\OOO$ over 
$(0,T)$.
\end{abstract} 
\baselineskip 18pt

\section{Introduction}

Let $\OOO \subset \R^d$, $d\in\mathbb{N}$, be a smooth bounded domain and 
$\nu = \nu(x)$ be the outward unit normal vector to $\ppp\OOO$ at 
$x \in \ppp\OOO$.  We set $\ppp_{\nu}v = \nabla v\cdot \nu$ and
$$
Q:= \OOO \times (0,T).
$$

For an arbitrarily given domain $D \subset \R^{d+1}$ in the $(x,t)$-space,
we set
$$
\Vert u\Vert_{H^{2,1}(D)}:= \left( \Vert u\Vert^2_{L^2(D)}
+ \Vert \nabla u\Vert^2_{L^2(D)} + \sumij \Vert \ppp_i\ppp_j u\Vert^2_{L^2(D)}
+ \Vert \ppp_tu\Vert^2_{L^2(D)}\right)^{\frac{1}{2}}.
$$
We consider a system of linearized mean field game equations:
$$
\left\{ \begin{array}{rl}
& \ppp_tu + \Delta u = Q_1u + Sv + F, \\
& \ppp_tv - \Delta v = Q_2(u,v) + \rho_0\Delta u + G \quad
\mbox{in $Q$,}
\end{array}\right.
                                \eqno{(1.1)}
$$
where $\rho_0$ is a constant (not necessarily positive) and 
$$
\left\{ \begin{array}{rl}
& \vert Q_1(u)(x,t)\vert \le C(\vert u(x,t)\vert + \vert \nabla u(x,t)\vert),\\
& \vert Q_2(u,v)(x,t)\vert \le C(\vert u(x,t)\vert + \vert v(x,t)\vert 
+ \vert \nabla u(x,t)\vert + \vert \nabla v(x,t)\vert).
\end{array}\right.
                                  \eqno{(1.2)}
$$
We assume 
$$
\vert S(v)(x,t) \vert \le C\vert v(x,t)\vert,  \quad (x,t) \in Q , \eqno{(1.3)}
$$
or
$$
\Vert Sv(\cdot,t)\Vert_{L^2(\OOO)} \le C\Vert v(\cdot,t)\Vert_{L^2(\OOO)},
\quad 0<t<T.                    \eqno{(1.4)}
$$
Moreover we let
$$
\ppp_{\nu}u = \ppp_{\nu}v = 0 \quad \mbox{on $\ppp\OOO\times (0,T)$}.
                                                                 \eqno{(1.5)}
$$

For small $\ep>0$, we set $Q_{\ep} := \OOO \times (\ep,\, T)$.

Let $u, v \in H^{2,1}(Q)$ satisfy (1.1)-(1.2) with (1.3) or (1.4).
Then
\\
{\bf Theorem 1.1.}
\\
{\it 
We assume 
$$
\Vert \nabla v(\cdot,0)\Vert_{L^(\OOO)}
\le M,
$$
where $M>0$ is a given constant.
For any $\ep \in (0,T)$, there exist constants $C=C(\ep,M) >0$ and 
$\theta=\theta(\ep, M) \in (0,1)$ 
such that 
$$
\Vert u\Vert_{H^{2,1}(\OOO \times (\ep,T))} 
+ \Vert v\Vert_{H^{2,1}(\OOO \times (\ep,T))} 
\le C(\Vert u(\cdot,T)\Vert_{H^1(\OOO)} + \Vert v(\cdot,T)\Vert_{L^2(\OOO)})
^{\theta}
$$
provided that $\Vert u(\cdot,T)\Vert_{H^1(\OOO)} 
+ \Vert v(\cdot,T)\Vert_{L^2(\OOO)}$ is sufficiently small.
}
\\
\vspace{0.2cm}

We refer to Lasry and Lions \cite{LL} 
for a convenient reference on the background of the mean field game equations 
(1.1)-(1.5) considered in this paper. 
As earlier works by estimates of Carleman's type, we refer to 
Klibanov \cite{Kl23}, Klibanov and Averboukh \cite{KlAv},
Klibanov, Li and Liu \cite{KLL1}, \cite{KLL2}.
The stability estimate in Theorem 1.1 was proved in 
Klibanov, Li and Liu \cite{KLL2}
for a nonlinear mean field game equations by a different
Carleman estimate which is not attached with second large parameter $\la>0$
such as in Theorem 2.1 stated below.  Thus the argument in \cite{KLL2}
requires extra estimation for $\rho_0\Delta u$ in the second equation in (1.1).
Thanks to the second large parameter $\la>0$, our proof is direct 
and is applicable to other cases 
where $\rho_0\Delta$ is replaced by arbitrary, not necessarily elliptic,
second-order partial differential operators. 
However, we here 
consider only a simple case (1.1) for describing our methodology.
On the other hand, the constant $C>0$ in the conclusion of Theorem 1.1 must be 
larger than \cite{KLL2},
Theorem 1.1 immediately implies the uniqueness for the backward problem 
for (1.1).  The theorem follows directly from the Carleman-type 
estimate (Theorem 2.1).  Our inverse problem is related to
backward parabolic equations in time, for which we can refer to many 
works and other methods may be available, 
and here we are limited only to two monographs: Ames and Straughan \cite{AS},
Payne \cite{Pa}.

In Klibanov and Averboukh \cite{KlAv}, the Lipschitz stability 
is proved in determining $u, v$ in $\OOO\times (0,T)$ by data 
$u(\cdot,T), v(\cdot,T), v(\cdot,0)$ in $\OOO$, while for the 
H\"older stability in $\OOO\times (\ep, T)$ with arbitrary 
$\ep > 0$ and the uniqueness in $\OOO \times (0,T)$, we do not need
data $v(\cdot,0)$ in $\OOO$.  Moreover Theorem 2.1 stated below 
readily implies such a global Lipschitz stability:
$$
\Vert u\Vert_{H^{2,1}(Q)} + \Vert v\Vert_{H^{2,1}(Q)}
\le C(\Vert u(\cdot,T)\Vert_{H^1(\OOO)} 
+ \Vert \nabla v(\cdot,0)\Vert_{L^2(\OOO)}
+ \Vert v(\cdot,T)\Vert_{L^2(\OOO)}),
$$
which can estimate the stronger norm $\Vert v\Vert_{H^{2,1}(Q)}$
than \cite{KlAv}.

As for a different kind of state determination, 
Klibanov, Li and Liu \cite{KLL1} proved conditional H\"older stability. 

We can choose the weight $\va(t) = e^{-\la t}$ with 
large constant $\la > 0$ in Theorem 2.1, and follow the proof of 
Theorem 1.1, so that we can obtain the H\"older stability 
in estimating $\Vert u\Vert_{H^{2,1}(\OOO \times (\ep,T))} 
+ \Vert v\Vert_{H^{2,1}(\OOO \times (\ep,T))}$ by  
$u(\cdot,0)$ and $v(\cdot,0)$ in $\OOO$.
In Klibanov \cite{Kl23}, the global Lipschitz stability 
by data $u(\cdot,0), u(\cdot,T), v(\cdot,0)$ in $\OOO$, but our result 
can prove the H\"older stability in determining $u,v$ in $\OOO\times
(\ep,T)$ only by $u(\cdot,0), v(\cdot,0)$ in $\OOO$.  Moreover 
our key Carleman estimate readily produces the global Lipschitz stability if
we can use $u(\cdot,0), u(\cdot,T), v(\cdot,0)$ in $\OOO$.
\section{Proof of Theorem 1.1}

We prove a key estimate of Carleman type with two large parameters.
We set 
$$
\va(t) = e^{\la t},
$$
where a constant $\la>0$ is chosen later.
\\
{\bf Theorem 2.1.}
\\
{\it 
There exists a constant $C>0$ such that 
\begin{align*}
& \int_Q \biggl\{ 
\vert \ppp_tu\vert^2 + \vert \Delta u\vert^2 + s\la\va\vert \nabla u\vert^2
+ s^2\la^2\va^2\vert u\vert^2\\
+ & \frac{1}{s\va}\left( \vert \ppp_tv\vert^2 + \vert \Delta v\vert^2
\right) + \la\vert \nabla v\vert^2 + s\la^2\va \vert v\vert^2
\biggr\} \weight dxdt\\
+ &\int_{\OOO} (s^2\la \vert u(x,0)\vert^2 
+ s\vert\nabla u(x,0)\vert^2)e^{2s} dx\\
+ &\int_{\OOO} (s\la \vert v(x,0)\vert^2 e^{2s}
+ \vert \nabla v(x,T)\vert^2 e^{2s\va(T)}) dx\\
\le & C\int_{\OOO} (s\la\va(T) \vert v(x,T)\vert^2 e^{2s\va(T)}
+ \vert \nabla v(x,0)\vert^2 e^{2s}) dx\\
+ & C\int_{\OOO} (s\la^2 \va(T)^2 \vert u(x,T)\vert^2 
+ s\va(T)\vert \nabla u(x,T)\vert^2) e^{2s\va(T)} dx
\end{align*}
for all large $s, \la > 0$.
}
\\

It is essential that the constant $\la>0$ in $e^{\la t}$ is activated as 
additional large parameter with $s>0$.
The weight $e^{s\va} = e^{se^{\la t}}$ is double exponential, so that 
the constant in the conclusion of Theorem 1.1 may be quite large, but the
proof is much easier.    
This weight was considered in Murray and Protter \cite{MP} for a different 
problem.

The proof follows from Theorem 9.1 in Yamamoto \cite{Y09}, and 
for convenience, I will give a sketch of the proof in Section 4.

Now we proceed to the proof of Theorem 1.1.
We fix $\la>0$ sufficiently large and we can delete the 
$\la$-dependence in Theorem 2.1:
\begin{align*}
& \int_Q \biggl\{ 
\vert \Delta u\vert^2 + \frac{1}{s}\vert \Delta v\vert^2
+ \vert \ppp_tu \vert^2 + \frac{1}{s}\vert \ppp_tv\vert^2\\
+ & (s\vert \nabla u\vert^2 + \vert \nabla v\vert^2
+ s^2\vert u\vert^2 + s\vert v\vert^2\biggr\} \weight dxdt\\
\le & C\int_{\OOO} e^{Cs}(s \vert v(x,T)\vert^2 
+ s\vert u(x,T)\vert^2 + \vert \nabla u(x,T)\vert^2) dx
\end{align*}
$$
+  C\int_{\OOO} \vert \nabla v(x,0)\vert^2 e^{2s} dx
                                                        \eqno{(2.1)}
$$
for all large $s > 0$.

We consider in $\OOO \times (\ep, T)$, and note 
$$
e^{2s\va(t)} \ge e^{2se^{\la\ep}} = e^{2s(\mu_0+1)},
$$
where we set 
$$
\mu_0:= e^{\la\ep} - 1 > 0,
$$
By $\Vert \nabla v(\cdot,0)\Vert_{L^2(\OOO)} \le M$, the estimate (2.1) implies
$$
e^{2s(\mu_0+1)} 
(\Vert u\Vert_{\HHHH(\OOO\times (\ep,T))}^2 
+ \Vert v\Vert_{\HHHH(\OOO\times (\ep,T))}^2)
\le Cse^{Cs}D^2 + CsM^2e^{2s}.
$$
Hence,
$$
\Vert u\Vert_{\HHHH(\OOO\times (\ep,T))}^2 
+ \Vert v\Vert_{\HHHH(\OOO\times (\ep,T))}^2
\le Cse^{Cs}D^2 + CM^2e^{-2s\mu_0}
\le Ce^{C_1s}D^2 + CM^2e^{-2s\mu_0}         \eqno{(2.2)}
$$
for all $s\ge s_0$: large constant.
Here and henceforth we set 
$$
D:= (\Vert u(\cdot,T)\Vert^2_{H^1(\OOO)} + \Vert v(\cdot,T)\Vert
_{L^2(\OOO)})^{\frac{1}{2}}.
$$
Replacing $s:= s+s_0$ and further choosing a large constant $C>0$, 
we have (2.2) for all $s>0$.
Without loss of generality, we can assume $D < M$.

Making the right-hand side of (2.2) large, as one possibility
we choose $s>0$ satisfying
$$
e^{C_1s}D^2 = e^{-2s\mu_0}M^2.
$$
Then, we have
$$
s = \frac{2}{C_1+2\mu_0}\log \frac{M}{D}.
$$
Then (2.2) yields the conclusion of Theorem 1.1 with
$$
\theta = \frac{2\mu_0}{C_1+2\mu_0} \in (0,1).
$$
Thus the proof of Theorem 1.1 is complete.
\section{Lipschitz stability in determination of state without
initial and final values}

We consider 
$$
\left\{ \begin{array}{rl}
& \ppp_tu + \Delta u = Q_1u + Sv + F, \\
& \ppp_tu - \Delta u = Q_2(u,v) + \rho_0\Delta u + G 
\quad \mbox{in $Q$}
\end{array}\right.
                                           \eqno{(3.1)}
$$
and
$$
\ppp_{\nu} u = \ppp_{\nu}v = 0 \quad \mbox{on 
$\ppp\OOO\times (0,T)$}.                        \eqno{(3.2)}
$$
Moreover we assume (1.2) for $Q_1$ and $Q_2$, and 
$$
\vert Sv(x,t)\vert \le C\vert v(x,t)\vert, \quad (x,t)\in Q.
$$
Let $\omega \subset \OOO$ be an arbitrarily chosen subdomain.

Then we can prove
\\
{\bf Theorem 3.1.}
\\
{\it
For arbitrarily given $\ep > 0$, we can find a constant $C_{\ep} >0$ such that 
$$
\Vert u\Vert_{\HHHH(\OOO\times (\ep,T-\ep))}
+ \Vert v\Vert_{\HHHH(\OOO\times (\ep,T-\ep))}
\le C_{\ep}(\Vert F\Vert_{L^2(Q)} + \Vert G\Vert_{L^2(Q)})
+ C_{\ep}(\Vert u\Vert_{L^2(\omega\times (0,T))}
+ \Vert v\Vert_{L^2(\omega\times (0,T))})
$$
and
$$
\Vert u(\cdot,t)\Vert_{L^2(\OOO)}
+ \Vert v(\cdot,t)\Vert_{L^2(\OOO)}
\le C_{\ep}(\Vert F\Vert_{L^2(Q)} + \Vert G\Vert_{L^2(Q)})
+ C_{\ep}(\Vert u\Vert_{L^2(\omega\times (0,T))}
+ \Vert v\Vert_{L^2(\omega\times (0,T))})
$$
for any $\ep \le t \le T-\ep$.
}
\\

This is the Lipschitz stability in the time interval
$(\ep,\, T-\ep)$.

{\bf Corollary.}
\\
{\it
In Theorem 3.1, we further assume that 
$$
Sv(x,t) = q(x,t)v(x,t), \quad (x,t) \in Q \quad 
\mbox{supp $q \cap \omega$ has interior points.}
$$
Then we can replace the right-hand side of the conlusions of 
Theorem 3.1 by 
$$
C(\Vert F\Vert_{L^2(Q)} + \Vert G\Vert_{L^2(Q)})
+ C\Vert u\Vert_{\HHHH(\omega\times (0,T))}.
$$
}

In the corollary, we note that
we do not need any extra data of $v$.
\\
\vspace{0.2cm}
\\
{\bf Proof.}
We set $t_0:= \frac{T}{2}$.

We know (e.g., Imanuvilov \cite{Ima}) that there exists $\eta \in 
C^2(\ooo{\OOO})$ such that 
$$
\eta\vert_{\ppp\OOO} = 0, \quad \eta > 0 \quad \mbox{in $\OOO$},
\quad \vert \nabla \eta \vert > 0 \quad \mbox{on 
$\ooo{\OOO \setminus \omega}$}.
$$
Set 
$$
\va(x,t) := \frac{e^{\la\eta(x)}}{t(T-t)}, \quad
\alpha(x,t) := \frac{e^{\la\eta(x)} - e^{2\la\Vert\eta\Vert_{C(\ooo{\OOO})}}}
{t(T-t)}.
$$
Then
\\
{\bf Lemma 3.1.}
\\
{\it
\begin{align*}
& \int_Q \biggl\{
\frac{1}{s\va}\left( \vert \ppp_tu\vert^2 
+ \sumij \vert \ppp_i\ppp_ju\vert^2\right)
+ s\va\vert \nabla u\vert^2 + s^3\va^3\vert u\vert^2
\biggr\} e^{2s\alpha} dxdt \\
\le& C\int_Q \vert (\ppp_t \pm \Delta)u\vert^2
e^{2s\alpha} dxdt 
+ C\int_{\omega \times (0,T)} s^3\va^3 \vert u\vert^2 e^{2s\alpha} dxdt
\end{align*}
for all $s \ge s_0$ and $u \in \HHHH(Q)$ satisfying 
$\ppp_{\nu}u = 0$ on $\ppp\OOO \times (0,T)$.
}
\\
{\bf Proof of Lemma 3.1.}

The case for $\ppp_t - \Delta$ is known; see e.g., 
Fursikov and Imanuvilov \cite{FI}, Imanuvilov \cite{Ima}.

We consider the case $\ppp_tu(x,t) + \Delta u(x,t) = F(x,t)$ for 
$x\in \OOO$ and $0<t<T$.
Setting $w(x,t):= u(x,T-t)$ for $0<t<T$.  Then
$\ppp_tw(x,t) = -(\ppp_tu)(x,T-t)$ and so
$$
\ppp_tw(x,t) - \Delta w(x,t)
= -(\ppp_tu)(x,T-t) - \Delta u(x,T-t)
= -F(x,T-t)
$$
and $\ppp_{\nu}w = 0$ on $\ppp\OOO\times (0,T)$.

Hence the Carleman estimate for $\ppp_t - \Delta$ yields
\begin{align*}
& \int_Q \biggl\{
\frac{1}{s\va}\left( \vert \ppp_tw(x,t)\vert^2 
+ \sumij \vert \ppp_i\ppp_jw\vert^2\right)
+ s\va\vert \nabla w\vert^2 + s^3\va^3\vert w\vert^2
\biggr\} e^{2s\alpha(x,t)} dxdt \\
\le& C\int_Q \vert F(x,T-t)\vert^2
e^{2s\alpha(x,t)} dxdt 
+ C\int_{\omega \times (0,T)} s^3\va^3 \vert w\vert^2 e^{2s\alpha(x,t)} dxdt,
\end{align*}
that is, using $\va(x,t) = \va(x,T-t)$ and
$\alpha(x,t) = \alpha(x,T-t)$, we obtain
\begin{align*}
& \int_Q \biggl\{
\frac{1}{s\va(x,T-t)}\left( \vert \ppp_tu(x,T-t)\vert^2 
+ \sumij \vert \ppp_i\ppp_ju(x,T-t)\vert^2\right)\\
+ & s2\va(x,T-t)\vert \nabla u(x,T-t)\vert^2 
+ s^3\va^3\vert u(x,T-t)\vert^2
\biggr\} e^{2s\alpha(x,T-t)} dxdt \\
\le& C\int_Q \vert F(x,T-t)\vert^2
e^{2s\alpha(x,T-t)} dxdt 
+ C\int_{\omega \times (0,T)} s^3\va^3 \vert u(x,T-t)\vert^2 
e^{2s\alpha(x,T-t)} dxdt.
\end{align*}
Changing the variables $t \mapsto \xi$ by $\xi:= T-t$, we can complete
the proof of Lemma 3.1.
$\blacksquare$
\\

Moreover 
\\
{\bf Lemma 3.2 (Imanuvilov and Yamamoto \cite{IY98}).}
\\
{\it
\begin{align*}
& \int_Q \biggl\{
\vert \ppp_tu\vert^2 
+ \sumij \vert \ppp_i\ppp_ju\vert^2
+ s^2\va^2\vert \nabla u\vert^2 + s^4\va^4\vert u\vert^2
\biggr\} e^{2s\alpha} dxdt \\
\le& C\int_Q s\va\vert (\ppp_t \pm \Delta)u\vert^2
e^{2s\alpha} dxdt 
+ C\int_{\omega \times (0,T)} s^4\va^4 \vert u\vert^2 e^{2s\alpha} dxdt
\end{align*}
for all large $s \ge s_0$ and $u \in \HHHH(Q)$ satisfying 
$\ppp_{\nu}u = 0$ on $\ppp\OOO \times (0,T)$.
}
\\

We apply Lemmata 3.2 and 3.1 to $u$ and $v$ respectively by replacing 
$F:= F + Q_1u + Sv$ and $G:= G + Q_2(u,v)+\rho_0\Delta u$:
\begin{align*}
& \int_Q \biggl\{
\vert \ppp_tu\vert^2 + \sumij \vert \ppp_i\ppp_ju\vert^2
+ s^2\va^2\vert \nabla u\vert^2 + s^4\va^4\vert u\vert^2
\biggr\} e^{2s\alpha} dxdt \\
\le& C\int_Q s\va \vert F\vert^2 e^{2s\alpha} dxdt 
+ C\int_Q s\va(\vert \nabla u\vert^2 + \vert u\vert^2 + \vert v\vert^2)
e^{2s\alpha} dxdt
\end{align*}
$$
+ C\int_{\omega \times (0,T)} s^4\va^4 \vert u\vert^2 e^{2s\alpha} dxdt
                                               \eqno{(3.3)}
$$
and
\begin{align*}
& \int_Q \biggl\{
\frac{1}{s\va}\left( \vert \ppp_tv\vert^2 
+ \sumij \vert \ppp_i\ppp_jv\vert^2\right)
+ s\va\vert \nabla v\vert^2 + s^3\va^3\vert v\vert^2
\biggr\} e^{2s\alpha} dxdt \\
\le& C\int_Q \vert G\vert^2 e^{2s\alpha} dxdt 
+ C\int_Q (\vert u\vert^2 + \vert \nabla u\vert^2
+ \vert v\vert^2 + \vert \nabla v\vert^2) e^{2s\alpha} dxdt
\end{align*}
$$
+ C\int_Q \vert \Delta u\vert^2 e^{2s\alpha} dxdt
+ C\int_{\omega \times (0,T)} s^3\va^3 \vert v\vert^2 e^{2s\alpha} dxdt
                                 \eqno{(3.4)}
$$
for all $s \ge s_0$.

Applying (3.3) to the third term on the right-hand side of (3.4), we obtain
\begin{align*}
& \int_Q \biggl\{
\frac{1}{s\va}\left( \vert \ppp_tv\vert^2 
+ \sumij \vert \ppp_i\ppp_jv\vert^2\right)
+ s\va\vert \nabla v\vert^2 + s^3\va^3\vert v\vert^2
\biggr\} e^{2s\alpha} dxdt \\
\le& C\int_Q (s\va\vert F\vert^2 
+ \vert G\vert^2) e^{2s\alpha} dxdt 
+ C\int_Q s\va(\vert u\vert^2 + \vert \nabla u\vert^2
+ \vert v\vert^2) e^{2s\alpha} dxdt\\
+ & C\int_{\omega \times (0,T)} s^4\va^4 \vert u\vert^2
e^{2s\alpha} dxdt\\
+& C\int_Q (\vert u\vert^2 + \vert \nabla u\vert^2
+ \vert v\vert^2 + \vert \nabla v\vert^2) e^{2s\alpha} dxdt\\
+ &C\int_{\omega \times (0,T)} s^3\va^3 \vert v\vert^2 e^{2s\alpha} dxdt.
\end{align*}
For large $s>0$, we can absorb the terms of lower-powers in $s$, so that 
we can obtain
\begin{align*}
& \int_Q \biggl\{
\frac{1}{s\va}\left( \vert \ppp_tv\vert^2 
+ \sumij \vert \ppp_i\ppp_jv\vert^2\right)
+ s\va\vert \nabla v\vert^2 + s^3\va^3\vert v\vert^2
\biggr\} e^{2s\alpha} dxdt \\
\le& C\int_Q (s\va\vert F\vert^2 
+ \vert G\vert^2) e^{2s\alpha} dxdt 
+ C\int_Q s\va(\vert u\vert^2 + \vert \nabla u\vert^2)
e^{2s\alpha} dxdt
\end{align*}
$$
+ C\int_{\omega \times (0,T)} (s^4\va^4 \vert u\vert^2
+ s^3\va^3\vert v\vert^2) e^{2s\alpha} dxdt               \eqno{(3.5)}
$$
for all large $s>0$.
Adding (3.3) and (3.5), and absorbing the terms of lower powers in $s\va$, 
we reach
\\
{\bf Theorem 3.2 (global Carleman estimate for the mean field game system).}
\\
{\it
\begin{align*}
& \int_Q \biggl\{ \vert \ppp_tu\vert^2
+ \sumij \vert \ppp_i\ppp_ju\vert^2 
+ s^2\va^2\vert \nabla u\vert^2 + s^4\va^4\vert u\vert^2\\
+& \frac{1}{s\va}\left( \vert \ppp_tv\vert^2 
+ \sumij \vert \ppp_i\ppp_jv\vert^2\right)
+ s\va\vert \nabla v\vert^2 + s^3\va^3\vert v\vert^2
\biggr\} e^{2s\alpha} dxdt \\
\le& C\int_Q (s\va\vert F\vert^2 
+ \vert G\vert^2) e^{2s\alpha} dxdt 
+ C\int_{\omega \times (0,T)} 
(s^4\va^4\vert u\vert^2 + s^3\va^3\vert v\vert^2)
   e^{2s\alpha} dxdt
\end{align*}
for all large $s>0$.
}
\\

Thanks to Theorem 3.2, we can readily complete the proof of Theorem 3.1.
Indeed, since 
$$
\alpha(x,t) \ge \frac{e^{\la\eta(x)} - e^{2\la\Vert \eta\Vert_{C(\ooo{\OOO})}} 
}{\ep (T-\ep)}
\ge \frac{-C_1}{\ep(T-\ep)} =: -C_2 < 0
$$
for all $x\in \ooo{\OOO}$ and $\ep \le t \le T-\ep$, we see that 
$$
e^{2s\alpha(x,t)} \ge e^{-2sC_2}, \quad x\in \ooo{\OOO}, \, 
\ep\le t \le T-\ep.
$$
Thus Theorem 3.2 completes the proof of Theorem 3.1.
$\blacksquare$
\section{Proof of Theorem 2.1.}

{\bf First Step.}

We consider 
$$
\left\{ \begin{array}{rl}
& \ppp_tv - \Delta v = \www{G} \quad \mbox{in $Q$},\\
& \ppp_{\nu}v = 0 \quad \mbox{on $\ppp\OOO\times (0,T)$}.
\end{array}\right.
                                      \eqno{(4.1)}
$$
Then we can prove
\begin{align*}
& \int_Q \biggl\{ \frac{1}{s\va}( \vert \ppp_tv\vert^2
+ \vert \Delta v\vert^2) + \la\vert \nabla v\vert^2
+ s\la^2\va\vert v\vert^2 \biggr\} \weight dxdt\\
+ & \int_{\OOO} (s\la \vert v(x,0)\vert^2 e^{2s}
+ \vert \nabla v(x,T)\vert^2 e^{2s\va(T)}) dx
\end{align*}
$$
\le C\int_{Q} \vert \www{G}\vert^2 \weight dxdt 
+ C\int_{\OOO} (s\la\va(T) \vert v(x,T)\vert^2 e^{2s\va(T)}
+ \vert \nabla v(x,0)\vert^2 e^{2s}) dx
                                              \eqno{(4.2)}
$$
for all large $s, \la > 1$ (Yamamoto \cite{Y09}).  
\\
{\bf Second Step.}
\\
Next we consider
$$
\left\{ \begin{array}{rl}
& \ppp_tu + \Delta u = \www{F} \quad \mbox{in $Q$},\\
& \ppp_{\nu}u = 0 \quad \mbox{on $\ppp\OOO\times (0,T)$}.
\end{array}\right.
                                      \eqno{(4.3)}
$$
Then we can prove
\begin{align*}
& \int_Q \biggl\{ \frac{1}{s\va}( \vert \ppp_tu\vert^2
+ \vert \Delta u\vert^2) + \la\vert \nabla u\vert^2
+ s\la^2\va\vert u\vert^2 \biggr\} \weight dxdt\\
+ & \int_{\OOO} (s\la \vert u(x,0)\vert^2 
+ \vert \nabla u(x,0)\vert^2) e^{2s} dx
\end{align*}
$$
\le C\int_{Q} \vert \www{F}\vert^2 \weight dxdt 
+ C\int_{\OOO} (s\la\va(T) \vert u(x,T)\vert^2 
+ \vert \nabla u(x,T)\vert^2) e^{2s\va(T)} dx
                                              \eqno{(4.4)}
$$
for all large $s, \la > 1$.  
\\
{\bf Proof of (4.4).}
The proof is similar to (4.2).  We set $w_1:= e^{s\va}u$ and 
$P_1w_1 = e^{s\va}(\ppp_t + \Delta)(e^{-s\va}w_1)$.
Then
$$
P_1w_1 = \ppp_tw_1 - s\la\va w_1 + \Delta w_1 = \www{F}e^{s\va}
$$ 
and so
\begin{align*}
& \Vert P_1w_1\Vert^2_{L^2(Q)}
\ge \int_Q \vert \ppp_tw_1\vert^2 dxdt 
+ 2\int_Q (\ppp_tw_1)\Delta w_1 dxdt
- \int_Q 2s\la\va w_1(\ppp_tw_1) dxdt \\
=: &\int_Q \vert \ppp_tw_1\vert^2 dxdt + \www{I_1} + \www{I_2}.
\end{align*}
Now we estimate
\begin{align*}
& \www{I_1} = -2\int_Q \nabla(\ppp_tw_1)\cdot \nabla w_1 dxdt
= -\int_Q \ppp_t(\vert \nabla w_1\vert^2) dxdt\\
=& \int_{\OOO} \left[ \vert \nabla w_1\vert^2\right]^0_T dx
= \int_{\OOO} (\vert \nabla w_1(x,0)\vert^2
- \vert \nabla w_1(x,T)\vert^2) dx
\end{align*}
and
\begin{align*}
& \www{I_2} = -\int_Q 2w_1(\ppp_tw_1)s\la\va dxdt
= - \int_Q \ppp_t(\vert w_1\vert^2)s\la\va dxdt \\
= & \int_{\OOO} \left[ s\la\va \vert w_1\vert^2\right]^0_T dx 
+ \int_Q s\la^2\va \vert w_1\vert^2 dxdt.
\end{align*}

Consequently,
$$
\int_Q \vert \ppp_tw_1\vert^2 dxdt 
+ \int_Q s\la^2\va \vert w_1\vert^2 dxdt
+ \int_{\OOO} (s\la\vert w_1(x,0)\vert^2 + \vert \nabla w_1(x,0)\vert^2)
dx
                                           \eqno{(4.5)}
$$
$$
\le C\int_Q \vert \www{F}\vert^2\weight dxdt 
+ C\int_{\OOO} (s\la\va(T)\vert w_1(x,T)\vert^2 + \vert \nabla w_1(x,T)\vert^2)
dx.
$$

Next we consider 
$$
w_1P_1w_1 = w_1\ppp_tw_1 - s\la\va \vert w_1\vert^2
+ w_1\Delta w_1   \quad \mbox{in $Q$}.
$$
Integrating over $Q$, we obtain
$$
  \int_Q \hhalf\ppp_t (\vert w_1\vert^2) dxdt 
- \int_Q s\la\va \vert w_1\vert^2 dxdt 
- \int_Q \vert \nabla w_1\vert^2 dxdt = \int_Q w_1(P_1w_1) dxdt.
$$
Hence
\begin{align*}
& -\hhalf \int_{\OOO} \left[ \vert w_1\vert^2\right]^T_0 dx 
+ \int_Q s\la\va \vert w_1\vert^2 dxdt 
+ \int_Q \vert \nabla w_1\vert^2 dxdt \\
=& -\int_Q w_1(P_1w_1) dxdt 
\le \int_Q \vert w_1\vert \vert P_1w_1\vert dxdt
dxdt,
\end{align*}
that is, the multiplication by $\la$ yields
\begin{align*}
& \la \int_Q \vert \nabla w_1\vert^2 dxdt 
+ \int_Q s\la^2\va \vert w_1\vert^2 dxdt 
+ \int_{\OOO} \la\left[ \vert w_1\vert^2\right]^0_T dx
\le & C \int_Q \la^2\vert w_1\vert^2 dxdt 
+ C\int_Q \vert \www{F}\vert^2 \weight dxdt.
\end{align*}
Here we used 
$$
\int_Q \la\vert w_1\vert \vert P_1w_1\vert dxdt
\le \hhalf\int_Q \la^2\vert w_1\vert^2 dxdt 
+ \hhalf\int_Q \vert P_1w_1\vert^2 dxdt.
$$
Adding (4.5), we see
\begin{align*}
& \int_Q \vert \ppp_tw_1\vert^2 dxdt 
+ \int_Q \left( 2s\la^2\va - C\la^2\right) \vert w_1\vert^2 dxdt
+ \int_Q \la \vert \nabla w_1\vert^2 dxdt\\
+ & \int_{\OOO} \{ (s\la+\la)\vert w_1(x,0)\vert^2 
+ \vert \nabla w_1(x,0)\vert^2 \} dx \\
\le &C\int_Q \vert \www{F}\vert^2\weight dxdt 
+ C\int_{\OOO} \{ (s\la\va(T)+\la)\vert w_1(x,T)\vert^2 
+ \vert \nabla w_1(x,T)\vert^2\} dx.
\end{align*}
Using 
$$
\frac{1}{s\va} \vert \ppp_tv\vert^2 e^{2s\va}
\le \frac{1}{s\va}\vert \ppp_tw\vert^2 + Cs\la^2\va\vert w\vert^2, 
$$
we can complete the proof of (4.4).
$\blacksquare$
\\
{\bf Third Step.}
\\
We set $u_1 := \va^{\hhalf}u$.  Then
$\ppp_tu_1 = \hhalf\la \va^{\hhalf}u + \va^{\hhalf}\ppp_tu$,
and so  
$$
\ppp_tu_1 + \Delta u_1 = \va^{\hhalf}\www{F} + \hhalf\la \va^{\hhalf}u.
                                             \eqno{(4.6)}
$$
Therefore the application of (4.4) to (4.6) yields
\begin{align*}
& \int_Q \biggl\{ \frac{1}{s\va}\va\left(
\left\vert \ppp_tu + \hhalf\la u\right\vert^2 + \vert \Delta u\vert^2
\right) + \la\va \vert \nabla u\vert^2 
+ s\la^2\va^2 \vert u\vert^2\biggr\} \weight dxdt \\
+ & \int_{\OOO} (s\la \vert u_1(x,0)\vert^2 + \vert \nabla u_1(x,0)\vert^2)
e^{2s} dx\\
= & \int_Q \biggl\{ \frac{1}{s\va}(\vert \ppp_tu_1\vert^2
+ \vert \Delta u_1\vert^2) + \la \vert \nabla u_1\vert^2
+ s\la^2\va\vert u_1\vert^2 \biggr\} \weight dxdt \\
+ & \int_{\OOO} (s\la \vert u_1(x,0)\vert^2 + \vert \nabla u_1(x,0)\vert^2)
e^{2s} dx\\
\le& C\int_Q \vert \www{F}\vert^2 \va \weight dxdt 
+ C\int_Q \la^2\va \vert u\vert^2 \weight dxdt
\end{align*}
$$
+ C\int_{\OOO} (s\la\va(T)\vert u_1(x,T)\vert^2 
+ \vert \nabla u_1(x,T)\vert^2)e^{2s\va(T)} dx.
                                              \eqno{(4.7)}
$$
Hence, choosing $s>0$ large, we can absorb the second term on the right-hand 
side of (4.7) into 
the left-hand side, we obtain
$$
\int_Q (\vert \ppp_tu\vert^2 - C\left\vert \hhalf\la u\right\vert^2
+ \vert \Delta u\vert^2 
+ s\la\va \vert \nabla u\vert^2 + s^2\la^2\va^2\vert u\vert^2) \weight dxdt 
                                         \eqno{(4.8)}
$$
\begin{align*}
+ &\int_{\OOO} (s^2\la\vert u(x,0)\vert^2 + s\vert \nabla u(x,0)\vert^2)
e^{2s} dx\\
\le & C\int_Q s\va\vert \www{F}\vert^2 \weight dxdt 
+ C\int_{\OOO} (s^2\la\va(T)^2\vert u(x,T)\vert^2 
+ s\va(T)\vert \nabla u(x,T)\vert^2)e^{2s\va(T)} dx\\
= & C\int_Q s\va\vert \www{F}\vert^2 \weight dxdt + CD(u(\cdot,T)).
\end{align*}
Here and henceforth, for short descriptions, we set 
$$
D(u(\cdot,T)):= \int_{\OOO} (s^2\la\va(T)^2\vert u(x,T)\vert^2 
+ s\va(T)\vert \nabla u(x,T)\vert^2)e^{2s\va(T)} dx
$$
with fixed large $\la, s>0$.
\\
{\bf Fourth Step.}
\\
We set $\www{F}:= Q_1u + Sv$ and $\www{G}:= Q_2(u,v) + \rho_0\Delta u$.
We note that 
\begin{align*}
& \vert \www{F}(x,t)\vert \le C(\vert Sv(x,t)\vert
+ \vert u(x,t)\vert + \vert \nabla u(x,t)\vert), \\
& \vert \www{G}(x,t)\vert \le C(\vert u(x,t)\vert 
+ \vert \nabla u(x,t)\vert + \vert v(x,t)\vert 
+ \vert \nabla v(x,t)\vert + \vert \Delta u(x,t)\vert)
\end{align*}
for $(x,t) \in Q$.  Moreover, in the case (1.4), we estimate  
$$
\int_Q s\va\vert Sv(x,t)\vert^2 \weight dxdt
= \int^T_0 s\va \weight \left( \int_{\OOO}
\vert Sv(x,t)\vert^2 dx \right) dt
= \int^T_0 s\va e^{2s\va} \Vert Sv(\cdot,t)\Vert_{L^2(\OOO)}^2 dt
$$
$$
\le C\int^T_0 s\va \weight \Vert v(\cdot,t)\Vert^2_{L^2(\OOO)}dt
= C\int_Q s\va \vert v\vert^2 \weight dxdt.    
$$
Therefore,
$$
\int_Q s\va\vert\www{F}\vert^2 \weight dxdt
\le C\int_Q s\va(\vert v\vert^2 + \vert u\vert^2
+ \vert \nabla u\vert^2) \weight dxdt.             \eqno{(4.9)}
$$
Moreover (4.8) implies
$$
 \int_Q \vert \www{G}\vert^2 \weight dxdt 
\le C\int_Q (\vert u\vert^2 
+ \vert \nabla u\vert^2 + \vert v\vert^2 
+ \vert \nabla v\vert^2 + \vert \Delta u\vert^2)\weight dxdt
$$
$$
\le  C\int_Q s\va\vert \www{F}\vert^2 \weight dxdt
+ C\int_Q (\vert v\vert^2 + \vert \nabla v\vert^2) \weight dxdt 
+ C\DDDDD.                                                \eqno{(4.10)}
$$
Adding (4.2) and (4.8), and applying (4.9) and (4.10), we obtain
\begin{align*}
& \int_Q \biggl\{ 
\vert \ppp_tu\vert^2 + \vert \Delta u\vert^2 
+ s\la\va \vert \nabla u\vert^2 + s^2\la^2\va^2\vert u\vert^2
+ \frac{1}{s\va}(\vert \ppp_tv\vert^2 + \vert \Delta v\vert^2)
+ \la\vert \nabla v\vert^2 + s\la^2\va\vert v\vert^2\biggr\} \weight dxdt \\
+& \int_{\OOO} (s\la\vert v(x,0)\vert^2 e^{2s}
+ \vert \nabla v(x,T)\vert^2e^{2s\va(T)}) dx\\
+ & \int_{\OOO} (s^2\la \vert u(x,0)\vert^2 
+ s\vert \nabla u(x,0)\vert^2)e^{2s} dx\\
\le& C\int_Q (s\va \vert \www{F}\vert^2 + \vert \www{G}\vert^2)
\weight dxdt 
+ C\int_{\OOO} (s\la\va(T)\vert v(x,T)\vert^2e^{2s\va(T)}
+ \vert \nabla v(x,0)\vert^2 e^{2s}) dx + CD(u(\cdot,T)) \\
\le & C\int_Q s\va(\vert v\vert^2 + \vert u\vert^2
+ \vert \nabla u\vert^2)\weight dxdt 
+ C\int_Q (\vert v\vert^2 + \vert \nabla v\vert^2) \weight dxdt \\
+& CD(u(\cdot,T)) 
+ C\int_{\OOO} (s\la\va(T)\vert v(x,T)\vert^2e^{2s\va(T)}
+ \vert \nabla v(x,0)\vert^2 e^{2s}) dx.
\end{align*}
Choosing $s, \la > 0$ large and noting that the powers of the terms
$\vert v\vert^2$ are $s\la^2\va$ on the left-hand side and 
$s\va$ and $1$ on the right-hand side, etc., we can absorb the first and 
the second terms on the right-hand 
side into the left-hand side,
we complete the proof of Theorem 2.1.
$\blacksquare$
\\
\vspace{0.2cm}
\\
{\bf Acknowledgements.}
The authors thanks Professor M.V. Klibanov (University of North Carolina at 
Charlotte) for valuable comments.
The work was supported by Grant-in-Aid for Scientific Research (A) 20H00117 
of Japan Society for the Promotion of Science.

\end{document}